\documentclass[12pt]{article}
\usepackage{latexsym}
\usepackage{amsfonts}
\usepackage{mathrsfs}
\usepackage{amssymb}
\usepackage{amsmath}
\usepackage{rawfonts}
\input{prepictex}
\input{pictex}
\input{postpictex}

\def\CP{{\mathbb C \mathbb P}}

\def\eea{\end{eqnarray*}}

\newtheorem{thm}{Theorem}
\newtheorem{prop}{Proposition}
\newtheorem{cor}{Corollary}
\newtheorem{lem}{Lemma}

\newenvironment{proof}{\medskip \noindent
{\bf Proof.}}{\hfill \rule{.5em}{1em}
\\}

\def\CC{\mathbb C}
\def\RR{\mathbb R}
\def\ZZ{\mathbb Z}

\def\+{\oplus}

\def\*{^{\ast}}

\def\RP{{\mathbb R \mathbb P}}

\newcounter{exmp}

\title{Twistors, Holomorphic Disks,\\ and Riemann Surfaces with Boundary}

\author{Claude LeBrun\thanks{Supported 
in part by  NSF grant DMS-0305865.}\\
Department  of Mathematics\\  SUNY, 
Stony Brook, NY 11794}

\date{}

\begin{document}
\maketitle

\begin{abstract} 
 Moduli spaces of  holomorphic disks in a complex
 manifold $Z$, with boundaries constrained to lie in a
 maximal  totally real submanifold $P$, have recently
 been found to underlie a number of geometrically
 rich twistor correspondences. The purpose of this  paper is to 
develop  a  general Fredholm regularity criterion for  
holomorphic curves-with-boundary $(\Sigma, \partial \Sigma )\subset (Z,P)$,
and  then show how this  applies, in particular,  to 
various moduli problems
 of twistor-theoretic interest. 
 \end{abstract}
 
\section{Introduction}
 
 Many interesting differential-geometric structures can best be 
 understood by means of {twistor correspondences}. 
 Here the main lesson is that moduli spaces of 
 compact complex curves $\Sigma$ in a complex manifold $Z$
 tend to carry  tautological  differential-geometric structures.
 Moreover,  the
 structures arising in this way   
  often actually represent the 
 general solution of some natural system of partial differential equations.

 The prototypical  construction  of this type   was first discovered by   
 Penrose, who called it the {\em nonlinear graviton}  \cite{pnlg}. 
 Suppose that $Z$ is a complex $3$-manifold, and let 
 ${\mathcal M}$ denote the moduli space of those  embedded 
 $\CP_1$'s  in  $Z$ 
 which have the same  normal bundle as does a
 projective line in $\CP_3$. 
Penrose discovered that  ${\mathcal M}$ is then a 
complex $4$-manifold, and naturally carries the holomorphic analog  of a 
conformal class of anti-self-dual metrics.

In order to get real geometries out of this fundamentally complex picture, 
one traditionally supposes that $Z$ is also equipped with an anti-holomorphic
involution $\sigma : Z\to Z$, and then focuses on the moduli space 
$M\subset {\mathcal M}$  of $\sigma$-invariant curves. 
The geometries that one obtains in this manner are necessarily 
real-analytic, but in many problems this  is  an expected feature of the general solution,  
anyway --- e.g.  
for reasons of  
elliptic regularity. For example, by considering complex 
$3$-folds $Z$ equipped with {\em free} anti-holomorphic involutions 
$\sigma$, Atiyah, Hitchin, and Singer 
\cite{ahs} showed that every self-dual Riemannian $4$-manifold
arises from the Penrose construction. Moreover, this conclusion holds not only 
 locally, but also  globally. 

Now one can also construct real-analytic  self-dual manifolds 
of metric signature $({+}{+}{-}{-})$ 
  by  instead equipping  a complex $3$-fold  $Z$  with an anti-holomorphic 
involution $\sigma$ with fixed-point set $P\neq  \varnothing$. 
However, the relevant differential-geometric problem now corresponds to an 
ultra-hyperbolic system of PDE rather than an elliptic one, and we
must therefore expect  most   solutions to 
 be of  low regularity.    Moreover, the focusing of  bicharacteristics 
 makes it exceedingly difficult to apply the  Penrose approach 
  to understand the global structure of solutions, even in the real-analytic context.

However, a new paradigm has recently emerged that 
 substantially sweeps  away  these obstacles. 
First,  forget about  the involution $\sigma$;   instead, focus on a
totally real submanifold $P\subset Z$,  
 previously thought to merely be the  fixed-point set of $\sigma$.
Second, forget about compact holomorphic curves; instead,  look for
holomorphic curves-with-boundary, where the boundary is constrained to 
lie on $P$.  Smooth solutions will arise from smooth $P$, while 
rougher solutions will arise from rougher $P$.

Lionel Mason and I  found this approach to be remarkably fruitful 
in our joint work on 
 Zoll surfaces \cite{lmzoll}
and 
split-signature $4$-manifolds \cite{lmfrei}. In this article, I will present a
new way of determining when a complex curve-with-boundary is stable
under small deformations of $(Z,P)$, and then show  how this tool can 
be used in the context of some interesting examples. My primary goal here is 
to indicate the wide applicability of these ideas to topics that now
seem ripe for further exploration. 

It is  perhaps worth
mentioning that quite different considerations have recently led 
 physicists    to  intensively 
 study  both closed \cite{witts} and open  \cite{berkwit} strings in twistor spaces.
 The distinction  between closed and open strings  precisely parallels
  the contrast between traditional twistor geometry  and the new kind 
 of twistor  correspondence explored here, with  
   the submanifold $P$ playing  the r\^ole of  an $m$-brane
  in  open string theory. I thus hope to convince you that  it's time we
   twistor geometers  opened  up and used our branes!

\section{Reflections on Kodaira}
 \label{reflex}

The systematic study of moduli spaces of complex submanifolds 
was largely instigated by Kodaira, who proved the following paradigmatic 
result \cite{kodsub}: 

\begin{thm}[Kodaira] \label{koda}
Suppose that  $X$ is a compact complex submanifold of a complex manifold
$Z$, and let $$N= [(T^{1,0}Z)|_X]/T^{1,0}X$$ be the normal bundle of 
$X$. If $H^1(X, {\mathcal O}(N))=0$, then the moduli space ${\mathcal M}$ of 
all compact complex submanifolds of $Z$   near 
$X$ 
is a complex manifold. Moreover, if $x\in  {\mathcal M}$
is the  the base-point representing  $X$, then there is a natural 
 isomorphism
   $T^{1,0}_x{\mathcal M}\cong H^0(X, {\mathcal O}(N))$. 
\end{thm}

In this picture, $X$ is assumed to be {\em embedded}
in $Z$, but a moment's thought immediately gives one a 
result in the non-embedded case, too:

\begin{cor}
Let $X$ be a compact complex manifold, and let 
$f: X\to Y$ be a holomorphic map. If 
$H^1(X, {\mathcal O} (f^*T^{1,0}Y))=0$, then there is a 
universal deformation for $f$ which is parameterized by 
a neighborhood of the origin in $H^0(X, {\mathcal O} (f^*T^{1,0}Y))$.
\end{cor}
\begin{proof}
Set $Z= X\times Y$, and embed $X$ in $Z$ as the graph of 
$f$. The normal bundle of this embedding is exactly 
$N= f^*T^{1,0}Y$. Now apply   Theorem \ref{koda}. 
\end{proof}

When $X=\CP_1$,
this corollary  is  frequently used in Mori theory \cite{mori}. 
 Indeed,   little  
harm is done  by even  invoking it  to  deform 
 {\em embedded} rational curves $\CP_1\subset Z$.
However, 
this corollary is ill  suited to the study  of embedded 
   complex curves $X$ of higher genus,    since 
in this setting   one often  has  $H^1(X, {\mathcal O}(T^{1,0}Z))\neq 0$
 even when   $H^1(X,{\mathcal O}(N))=0$. 

Many applications of   twistor ideas depend on   
 the persistence of families of complex submanifolds after deformation
of the ambient complex manifold. Fortunately, 
a  minor modification
 \cite{kodstab} of  
Theorem \ref{koda} provides a criterion for guaranteeing 
 the survival of complex submanifolds
 in this context: 

\begin{thm}[Kodaira] \label{coda}
Suppose that  $X\subset Z$ is a compact complex submanifold 
whose   normal bundle satisfies  $H^1(X, {\mathcal O}(N))=0$.
Then any small deformation $Z^\prime$ of $Z$  contains
an   $h^0(X, {\mathcal O}(N))$-complex-dimensional  family 
of compact complex submanifolds, obtained by  deforming $X$. 
\end{thm}

Now these theorems do not depend  at all on the dimension of 
$X$. Nonetheless, the special case in which $X$ is a Riemann surface 
enjoys  a  somewhat privileged status. 
For example,  the same framework then  works even for {\em pseudo-holomorphic} curves 
 in {\em almost}-complex manifolds $Z$, and  
 Gromov \cite{gromsym} was able to make systematic
use of this observation  to prove    a   family of truly revolutionary 
results on the structure of  symplectic manifolds. 
However, it is yet 
{\em another}  special feature of the  $1$-dimensional case which will
concern us here; namely, as will be explained below, 
  these ideas can also  be naturally generalized so as to 
handle Riemann surfaces with non-empty boundary \cite{forst,glob,mcsalt,ohrk}.

 Let $Z$ be a complex $n$-manifold, and let $J$ denote its complex structure tensor. 
 Suppose that 
   $P$ is  a differentiable 
submanifold of $Z$  of real dimension $n$. We will then say that
$P$ is a  (maximal) totally real submanifold if $T_yP\cap J(T_yP) =0$ for
all $y\in P$, so   that $TZ|_P= TP\oplus J(TP)$. 

Now suppose that $\Sigma$ is a compact complex curve-with-boundary,
and that $\Sigma \hookrightarrow Z$ is a holomorphic embedding 
that is differentiable up to the boundary. Also assume that 
$\partial \Sigma \subset P$, where $P$ is a maximal totally real submanifold
of differentiability class $C^{k+4}$, $k\geq 1$. A regularity theorem of
Chirka \cite{chirka} then asserts that $\Sigma \hookrightarrow Z$ 
is actually $C^{k+3}$. Our goal here is 
to understand the space of nearby $C^1$ holomorphic curves $\Sigma^\prime$
with $\partial \Sigma^\prime \subset P$. However,   Chirka's regularity result  tells us that 
this
is  equivalent,  for example, to studying holomorphic curves  of class  $C^{k+1, \alpha}$
for any chosen $\alpha\in (0,1)$. Needless to say, we could now elect to set 
 $k=1$, once and for all, 
but I will leave  the choice of $k$ up to the reader, as doing so 
may actually  clarify certain aspects of the argument. 

Let us next choose a $C^{k+3}$ open surface ${\mathcal S}\subset Z$
containing $\Sigma$
and a $C^{k+3}$ Riemannian metric $g$ on $Z$ with respect to which 
$P$ is totally geodesic. Let $E\subset TP|_{\partial \Sigma}$ be the orthogonal
complement of $T\partial \Sigma$ relative to $P$, and, after shrinking $\mathcal S$ if necessary, 
 let 
$\tilde{N}\subset T^{1,0}Z|_{\mathcal S}$ be 
  a $C^{k+2}$  complex sub-bundle whose
  real part $\Re e\tilde{N}$ is complementary to $T{\mathcal S}$ and agrees with 
  $E\oplus J(E)$ along $\partial \Sigma$. Applying the geodesic spray of 
  $g$ to $\Re e\tilde{N}$ and invoking the inverse function theorem, 
  we thus obtain a  $C^{k+2}$ diffeomorphism
$\Phi$ 
between some neighborhood ${\mathcal V}\subset N$ of $\Sigma \subset
0_{\mathcal S}$ and an 
open subset $\mathcal U\subset Z$. Notice,  moreover, that we have arranged that  
$$
P\cap {\mathcal U} = \Phi (E\cap {\mathcal V}). 
$$

Now  notice that $N=\tilde{N}|_\Sigma$ can canonically be 
identified with the normal bundle $T^{1,0}Z/T^{1,0}\Sigma$ of $\Sigma$.
Choose some  inner product and connection on $N$, and 
let $C^{k+1,\alpha} (N, E) $
denote the Banach space of $C^{k+1,\alpha}$ sections of
$N$ whose boundary values are sections of 
$E\to \partial \Sigma$; let 
$C^{k+1,\alpha} (N, E)_\varepsilon \subset C^{k+1,\alpha} (N, E) $
be the $\varepsilon$ ball about $0$ in this Banach space. If  
$\varepsilon$ is sufficiently small, the graph of any $f\in C^{k+1,\alpha} (N, E)_\varepsilon$  is 
contained in 
$\mathcal V$, and so is sent by $\Phi$ to a $C^{k+1,\alpha}$ surface
$\Sigma^\prime \subset Z$ with  $\partial \Sigma^\prime \subset 
P$. Now let ${\mathbf V}^{1,0}\subset T_\CC N$ denote the vertical vectors of type 
$(1,0)$, and notice that this is naturally isomorphic to the pull-back 
of $N$ to its own total space. By possibly shrinking
our neighborhood $\mathcal V$ of $\Sigma$, we then 
have 
$$
T^{0,1}Z\cap \Phi_* {\mathbf V}^{1,0} = 0. 
$$
On the other hand, by shrinking $\varepsilon$ if necessary, 
we may also arrange that  
$$
T^{0,1}Z\cap (\Phi\circ f)_*[T^{1,0}\Sigma]=0
$$ 
 for all $f\in C^{k+1,\alpha} (N, E)_\varepsilon$. 
 Hence 
  $$T_\CC Z= T^{0,1}Z + \Phi_*{\mathbf V}^{1,0}+  (\Phi\circ f)_*[T^{1,0}\Sigma]$$
 at each  point of the image of any $f\in C^{k+1,\alpha} (N, E)_\varepsilon$.
 Composing $(\Phi \circ f)_*$ with the  projection 
 $T_\CC Z\to  \Phi_*{\mathbf V}^{1,0}$  thus defines a linear map 
 $T_\CC \Sigma \to N$ for every such $f$; and since this linear map 
 kills $T^{1,0}\Sigma$ by construction, it may be viewed as 
 a $(0,1)$-form ${\mathscr D}f$  with values in $N$. 
It is now easy to see that 
 $${\mathscr D}: C^{k+1,\alpha} (N, E)_\varepsilon \to C^{k,\alpha} (\Lambda^{0,1}\otimes N)$$
 is a differentiable  map 
 of Banach manifolds, and that the linearization of ${\mathscr D}$ at
 $0$ is exactly the canonical  operator 
$$\overline{\partial} :   C^{k+1,\alpha} (N, E) \to C^{k,\alpha} (\Lambda^{0,1}\otimes N),$$
obtained by remembering that   $N=T^{1,0}Z/T^{1,0}\Sigma$ is  a holomorphic
vector bundle over  $\Sigma$. 
Notice, moreover,  that ${\mathscr D}^{-1}(0)$ exactly consists of those holomorphic curves
$(\Sigma^\prime , \partial \Sigma^\prime )\hookrightarrow (Z,P)$
which are  sufficiently near $\Sigma$. 

% 
% Thus the projection $(\Phi\circ f)_*^{1,0}[T^{1,0}\Sigma]$ of
% $(\Phi\circ f)_*[T^{1,0}\Sigma]$ into $T^{1,0}Z= T_\CC Z/T^{0,1}Z$
% is a line bundle on $\Sigma$ for each such $f$, and the pull-back of 
% $$T^{1,0}Z/(\Phi\circ f)_*^{1,0}(T^{1,0}\Sigma)
% =T_\CC Z/\left[ T^{0,1}Z+(\Phi\circ f)_*^{1,0}(T^{1,0}\Sigma)\right]
% $$
% is therefore a complex vector bundle  of rank $(n-1)$ on $\Sigma$ for every
% $f\in C^{k+1,\alpha} (N, E)_\varepsilon$. What is more, we can actually
% identify each  such vector bundle with the  normal bundle 
% $N=T^{1,0}Z/T^{1,0}\Sigma$ of $\Sigma$,; namely,  we can  identify the pull-back
% of $N$ to $N$ with the bundle $V^{1,0}\subset T_\CC N$ consisting of vertical vectors of type $(1,0)$,
% and, by shrinking $\varepsilon$ if necessary,  
% we have 

I now want to explain a simple geometric trick  that not only proves that
this linearized operator is Fredholm, but actually provides a practical method
of precisely calculating the its kernel and cokernel.  The key idea is to
first construct the {\em abstract double} of our Riemann surface. 
That is, we begin with our Riemann-surface-with-boundary $\Sigma$

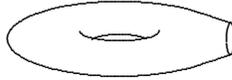
\begin{center}
\mbox{
\beginpicture
\setplotarea x from 0 to 283, y from 0 to 60
\ellipticalarc axes ratio 3:1  270 degrees from 145 40
center at 115 30
%\ellipticalarc axes ratio 3:1  -270 degrees from 180 40
%center at 210 30
\ellipticalarc axes ratio 4:1 -180 degrees from 130 33
center at 115 33
\ellipticalarc axes ratio 4:1 145 degrees from 125 30
center at 115 29
%\ellipticalarc axes ratio 4:1 180 degrees from 195 33
%center at 210 33
%\ellipticalarc axes ratio 4:1 -145 degrees from 200 30
%center at 210 29
\ellipticalarc axes ratio 1:4 360 degrees from 157 36
center at 157 30
%\ellipticalarc axes ratio 1:3 180 degrees from 168 36
%center at 168 30
{\setlinear 
\plot 145 40        157 36   /
\plot 145 20        157 24  /
%\plot 180 40        168 36   /
%\plot 180 20        168 24   /
}
\endpicture
}
\end{center}
and then attach a mirror-image copy $\overline{\Sigma}$ to $\Sigma$ 
along $\partial \Sigma$:

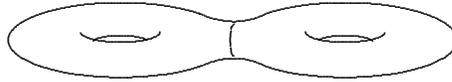
\begin{center}
\mbox{
\beginpicture
\setplotarea x from 0 to 290, y from 0 to 60
\ellipticalarc axes ratio 3:1  270 degrees from 150 40
center at 120 30
\ellipticalarc axes ratio 3:1  -270 degrees from 175 40
center at 205 30
\ellipticalarc axes ratio 4:1 -180 degrees from 135 33
center at 120 33
\ellipticalarc axes ratio 4:1 145 degrees from 130 30
center at 120 29
\ellipticalarc axes ratio 4:1 180 degrees from 190 33
center at 205 33
\ellipticalarc axes ratio 4:1 -145 degrees from 195 30
center at 205 29
\ellipticalarc axes ratio 1:4 180 degrees from 163 36
center at 163 30
{\setquadratic 
\plot 150 40    163 37    175 40   /
\plot 150 20    163 23    175 20  /
}
\endpicture
}
\end{center}
Let ${\mathbb X}= \Sigma \cup_{\partial \Sigma} \overline{\Sigma}$ denote
this double, and notice that it comes equipped with an
anti-holomorphic involution 
$$\rho : {\mathbb X}\stackrel{\overline{\mathcal O}}{\longrightarrow} {\mathbb X},
~~~~~\rho^2= \mbox{id}_{{\mathbb X}},$$
obtained by interchanging $\Sigma$ and $\overline{\Sigma}$.
This is true  because 
 $\overline{\Sigma}$ is by definition simply  
$\Sigma$, equipped with its {\em conjugate} complex structure.  

Now notice that  $\overline{N}$   is a 
a holomorphic vector bundle on $\overline{\Sigma}$. 
On the other hand, both $N$ and $\overline{N}$ restrict
to $\partial \Sigma$ as complexifications of the real
vector bundle $E= TP/T(\partial \Sigma )$. 
We can therefore  construct a 
complex vector bundle ${\mathscr N}\to {\mathbb X}$ 
by attaching  $\overline{N}\to \overline{\Sigma}$ to 
$N\to \Sigma$ in such a manner that $E$ is sent
to itself by the identity: 
$$
\begin{array}{ccccc}
    {\mathscr N} &  = &  N &\cup_{E \otimes \CC}& \overline{N}  \\
     \downarrow &&\downarrow && \downarrow \\
      {\mathbb X}&=& \Sigma&\cup_{\partial \Sigma}& \overline{\Sigma}
\end{array}
$$
We can make this into a holomorphic vector bundle by 
taking the complex structure tensor on its total space to be that of $N$ over  
$\Sigma$ and that of $\overline{N}$ over $\overline{\Sigma}$. 
Of course, we still need to check that  this gives us a locally 
trivial structure 
 in the vicinity of   $\partial \Sigma$. To see this, first recall 
  that we  arranged for
 $\Sigma \hookrightarrow Z$ to be  at least 
   $C^2$ up to the boundary, so that, even near  $\partial \Sigma$,  the bundle
   $N$ has $C^1$ local holomorphic
   trivializations  induced by local holomorphic trivializations of 
   $T^{1,0}Z$. The integrable almost-complex structure on 
   total space of $N$  is therefore $C^1$ up to the boundary. By reflection, the 
   conjugate complex structure of $\overline{N}$ is therefore $C^1$ up
   to the boundary, too. Now,  the manner in which we glue $N$ and 
   $\overline{N}$ together to make a $C^1$ manifold is 
   exactly chosen so that these two almost complex structures
   agree along the interface. Hence the total space of ${\mathscr N}$
  carries an induced almost-complex structure which is at least Lipschitz. 
   However, a result of 
Nijenhuis and Woolf \cite{nijac} asserts that a Lipschitz   almost-complex manifold  
contains 
 pseudo-holomorphic curves tangent to any given complex tangent line in its
 tangent space.  
But in our case  the generic such pseudo-holomorphic curve is necessarily  the graph
of a local  section of $\mathscr N$, and we thus obtain  enough
local holomorphic sections 
 of $\mathscr N$ to generate local holomorphic
trivializations, even  near points of $\partial \Sigma$. 
Thus ${\mathscr N}\to {\mathbb X}$ really is a holomorphic vector
bundle, as claimed.
%\footnote{This could of course also be shown 
%via    the ``nuclear option''  of 
% taking $\alpha > 1/2$ and  appealing to  the 
% $C^{0,\alpha}$  Newlander-Nirenberg theorem 
%of  Hill and Taylor \cite{hiltay}.}

 Now notice that, by construction,  the total 
space of ${\mathscr N}$ carries a tautological 
involution $\varrho$ which covers $\rho$:
\begin{eqnarray*}
 {\mathscr N}  & \stackrel{\varrho}{\longrightarrow}  &  {\mathscr N}  \\
 \downarrow &  & \downarrow \\
 {\mathbb X}& \stackrel{\rho}{\longrightarrow}  & {\mathbb X} 
\end{eqnarray*}
The fixed-point set of this involution is exactly the given sub-bundle
$E\subset {\mathscr N}|_{\partial \Sigma}$. A somewhat surprising 
consequence of this is that, no matter how rough $E$ may have appeared in
our original picture, it is actually real-analytic as a subspace of $\mathscr N$.
Indeed, we can even find holomorphic local trivializations of
$\mathscr N$ near any point of $\partial \Sigma$ in which 
$E$ becomes the trivial bundle with fiber $\RR^n$ as a 
sub-bundle of the trivial bundle with fiber $\CC^n$. 
Such special local trivializations will play a prominent r\^ole in what follows. 

For this reason,  it is important that we now check that these special 
 local trivializations are actually  of class $C^{k+1,\alpha}$
relative to  the na\"{\i}ve local trivializations of $N$. To see this, 
let $\{h_j\}$ be a local holomorphic frame
for ${\mathscr N}$ whose real span along $\partial \Sigma$ is $E$. 
Working instead with respect to a local trivialization of $N\to \Sigma$ induced by 
 some local holomorphic vector fields on $Z$, the assumed regularity of $P$ allow us to 
choose a $C^{k+2}$ local frame $\{ e_k\}$ for $N\to \Sigma$ whose 
real span  along $\partial \Sigma$ is $E$. Let $\psi$ be a
smooth bump function supported in the common domain of $\{e_j\}$ and
$\{ h_k\}$. 
We then have 
$$\psi e_j = \sum_k c_{jk}h_k$$
where $c_{jk}$ is real along $\partial \Sigma$, and where 
$\overline{\partial} c_{jk}$ is $C^{k+1}$. By taking local coordinates
on $\Sigma$, we can then view each $c_{jk}$ as a compactly supported
function on the upper half-plane, and then convert this into a 
smooth function on the $2$-disk $D^2$ by a applying a M\"obius
transformation. Expressing 
$c_{jk} = a_{jk} + i b_{jk}$
in terms of its real and imaginary parts,  $b_{jk}$ then vanishes
along the $\partial D^2$, while $\Delta b_{jk}$ is of class 
$C^{k-1}$ on $D^2$. Elliptic regularity for the Dirichlet problem   \cite{giltrud} thus 
predicts that $b_{jk}$ is of class $C^{k+1,\alpha}$ for any $\alpha\in (0,1)$,
and the fact that $ da_{jk}+ J(d b_{jk})$ is $C^{k+1}$ then implies that the 
$a_{jk}$ must be of class $C^{k+1,\alpha}$, too. It follows that the 
$\{ h_{j}\}$ are also $C^{k+1,\alpha}$,  as claimed.

Next, notice that  $\varrho$  induces complex-anti-linear involutions
$$\varrho^*: H^j ({\mathbb X}, {\mathcal O} ({\mathscr N}))\to 
H^j ({\mathbb X}, {\mathcal O} ({\mathscr N})),~~~ j =0,1.$$ 
Let $H^j_\varrho ({\mathbb X}, {\mathcal O} ({\mathscr N}))$
denote the $(+1)$-real-eigenspace of $\varrho^*$, so that 
$$H^j ({\mathbb X}, {\mathcal O} ({\mathscr N})) =H^j_\varrho ({\mathbb X}, {\mathcal O} ({\mathscr N}))\oplus iH^j_\varrho ({\mathbb X}, {\mathcal O} ({\mathscr N}))$$
as a real vector space. With this notation in hand, we can now  formulate our key technical
result:

\begin{lem}
\label{lift}
 For any integer $k\geq 1$ and any $\alpha \in (0,1)$,  
the linear operator
$$\overline{\partial} :   C^{k+1,\alpha}(\Sigma ; N,E) \to C^{k,\alpha}(\Sigma ;\Lambda^{0,1}\otimes N)$$
is Fredholm, with kernel canonically isomorphic to 
$H^0_\varrho({\mathbb X}, {\mathcal O} ({\mathscr N}))$ and cokernel canonically isomorphic to 
$H^1_\varrho({\mathbb X}, {\mathcal O} ({\mathscr N}))$. 
In particular, $\ker \overline{\partial}$ has  real dimension 
$h^0({\mathbb X}, {\mathcal O} ({\mathscr N}))$, while 
$\mbox{\rm coker } \overline{\partial}$ 
has  real dimension
$h^1({\mathbb X}, {\mathcal O} ({\mathscr N}))$. 
\end{lem}

\begin{proof}
First, let us compute the kernel of $\overline{\partial}$. If $f\in C^{k+1,\alpha}(\Sigma ; N,E)$, 
define a continuous section $\mathfrak f$ of $\mathscr N \to \mathbb X$ by 
$$
\mathfrak f = \left\{ \begin{array}{cc}
    f  &    \mbox{on } \Sigma\\
     \varrho^*{f} &   \mbox{on } \overline{\Sigma}
\end{array}\right.
$$
This is  well defined and continuous because, by assumption, 
$f=\varrho^*{f}$ along $\partial \Sigma$. Now if  $f\in \ker \overline{\partial}$,
$\mathfrak f$ is continuous up to $\partial \Sigma$ and holomorphic
on its complement, and so is holomorphic on all of $\mathbb X$
by the reflection principle.
Moreover, $\mathfrak f$ is invariant under the action of $\varrho^*$ by 
construction. Hence $\mathfrak f \in H^0_\varrho({\mathbb X}, {\mathcal O} ({\mathscr N}))$.
Since any element of $H^0_\varrho({\mathbb X}, {\mathcal O} ({\mathscr N}))$
conversely restricts to $\Sigma$ as 
an element of 
$C^{k+1,\alpha}(\Sigma ; N,E)$ 
which is killed by $\overline{\partial}$, we thus conclude that 
$H^0_\varrho({\mathbb X}, {\mathcal O} ({\mathscr N}))$ can naturally be 
identified with the kernel of the  operator.

Now, what is the   image of the 
operator? Given $\phi\in C^{k,\alpha}(\Sigma ;\Lambda^{0,1}\otimes N)$, define
a section of $\Lambda^{0,1}\otimes {\mathscr N}\to
{\mathbb X}$ by 
$$
\varphi = \left\{ \begin{array}{cc}
    \phi  &    \mbox{on } \Sigma\\
     \varrho^*{\phi} &   \mbox{on } {\mathbb X}-{\Sigma}
\end{array}\right.
$$
Of course, this $\varphi$ may not be continuous, but at any rate it is certainly $L^\infty$, and 
in particular may be considered as 
a twisted distribution-valued $(0,1)$-form.   Since Dolbeault cohomology 
can  be computed using currents \cite{gunningrs} {\em en lieu}
 of smooth forms, it follows that there is  a well defined   cohomology class
$[\varphi]\in H^1 ({\mathbb X},  {\mathcal O} ({\mathscr N}))$ which precisely 
measures the 
obstruction to writing $\varphi$  as 
$$
\varphi = \overline{\partial} {\mathfrak f}
$$
for some 
 distributional section $\mathfrak f$ of ${\mathscr N}\to {\mathbb X}$; moreover, 
 $[\varphi]\in H^1_\varrho ({\mathbb X},  {\mathcal O} ({\mathscr N}))$, since, 
  by construction,  $\varphi$ is $\varrho$-invariant almost everywhere.  
 We thus have a continuous linear map  
\begin{eqnarray*}
\Pi : C^{k,\alpha}(\Sigma ;\Lambda^{0,1}\otimes N) & \to  &  
H^1_\varrho ({\mathbb X},  {\mathcal O} ({\mathscr N}))\\
\phi & \mapsto & [\varphi ] 
\end{eqnarray*}
Moreover, this map is a surjection, 
 since $\mathbb X$ has a Stein cover consisting of any small  neighborhood $U$ of 
 $\Sigma$ and its conjugate $\overline{U}$,  and every element of 
 $H^1_\varrho ({\mathbb X},  {\mathcal O} ({\mathscr N}))$
 can be expressed as $[\varphi ]$ for a smooth section $\phi$ of 
 $\Lambda^{0,1}\otimes N\to \Sigma$ obtained by cutting off a 
  \v{C}ech representative $\in \Gamma (U\cap \overline{U} , {\mathcal O}({\mathscr N}))$
 with a bump function  and then restricting  to $\Sigma$.

Now notice that  $\overline{\partial} C^{k+1,\alpha}(\Sigma ; N,E)\subset \ker \Pi$. Indeed, 
 if $\phi = \overline{\partial} f$ for some
 $f \in C^{k+1,\alpha}(\Sigma ; N,E)$, the continuous section 
 $$
\mathfrak f = \left\{ \begin{array}{cc}
    f  &    \mbox{on } \Sigma\\
     \varrho^*{f} &   \mbox{on } \overline{\Sigma}
\end{array}\right.
$$ 
of $\mathscr N$ then satisfies $\overline{\partial} {\mathfrak f}=\varphi$ in the
distributional sense, and  $[\varphi]= \Pi (\phi )$ therefore vanishes. 

To finish the proof, it therefore suffices to show that 
$\ker \Pi \subset \overline{\partial} C^{k+1,\alpha}(\Sigma ; N,E)$. 
Thus,  suppose that we are given some  
$\phi \in C^{k,\alpha}(\Sigma ;\Lambda^{0,1}\otimes N)$ for which 
 $[\varphi] =0 \in H^1 ({\mathbb X},  {\mathcal O} ({\mathscr N}))$. It then follows that 
$\varphi = \overline{\partial} {\mathfrak f}_0$ for
some distributional ${\mathfrak f}_0$. 
Since $\varphi$ is $L^p$ for any $p$, 
elliptic regularity then tells us that 
${\mathfrak f}_0 \in L_1^p({\mathbb X}, {\mathscr N})$ for any $p$. 
Taking $p> 2$,  thus have  
${\mathfrak f}_0\in C^{0}({\mathbb X}, {\mathscr N})$
by the Sobolev embedding theorem.
 Setting  $\mathfrak f = ({\mathfrak f}_0+ \varrho^*{\mathfrak f}_0)/2$
we then have $\varphi = \overline{\partial}{\mathfrak f}$, where  $\mathfrak f$
is a $\varrho^*$-invariant continuous section of $\mathscr N$,  and so takes values 
in $E$ along $\partial \Sigma$ by
$\varrho^*$ invariance. Letting $f$ denote ${\mathfrak f}|_\Sigma$, we then
have $f\in C^{0}(\Sigma ; N, E)$,  and the point worth emphasizing is that 
$f|_{\partial \Sigma}$  is a section of $E$. Moreover, Schauder theory \cite{giltrud} tells us that 
$f$ is of class $C^{k+1, \alpha}$
on the interior of $\Sigma$.
It therefore only remains to show that $f$ is $C^{k+1,\alpha}$ in
the vicinity of any boundary point. 

Since this last  issue is completely local, we can multiply $f$ by a smooth bump function
supported in the domain of a special local trivialization near a given boundary point,
and so 
  obtain a weak solution of the 
equation 
$$\overline{\partial} \hat{f} = \hat{\phi}$$
where  $\hat{f}$ is  a compactly supported ${\mathbb C}^n$-valued continuous function
on the upper half-plane which is $\RR^n$-valued along the real axis, and where
$\hat{\phi}$ is a ${\mathbb C}^n$-valued $(0,1)$-form of class $C^{k,\alpha}$.
We now identify the upper half-plane with the $2$-disk $D$ via a M\"obius transformation. 
Expressing the $j^{\rm th}$ component of $\hat{f}$ in terms of its real and imaginary parts
$$\hat{f}_j = u + iv,$$ and correspondingly expressing the $j^{\rm th}$ component of $\hat{\phi}$
as $$\hat{\phi}_j= \frac{\alpha + i\beta}{2} d\overline{z},$$ it then follows that $v$
is a weak solution of the Dirichlet problem 
$$\Delta v = h~~
\mbox{   on }D,~~~~v=0 ~~\mbox{   on }\partial D, $$
where 
$$
h =  \frac{\partial \alpha}{\partial y} - \frac{\partial \beta}{\partial x}
$$
is of class $C^{k-1,\alpha}$. 
It follows \cite{giltrud} that  $v$ is of class $C^{k+1,\alpha}$, as desired. But we also have
$$du = -Jdv + \alpha ~dx  + \beta ~dy$$
so this implies  that $u$ is of class $C^{k+1,\alpha}$, too. Hence
$f$ is everywhere of class  $C^{k+1,\alpha}$, as claimed .
\end{proof}

Notice that the same reasoning actually applies to any holomorphic
vector bundle on any $\Sigma$ and any maximal real sub-bundle of its
restriction to $\partial \Sigma$. 
Of course,  the precise computations of the kernel and cokernel 
are delicate in nature, but the fact that the map is Fredholm is stable under
perturbation by compact operators. Thus the Fredholm property holds for quite  general 
first-order operators of Cauchy-Riemann type with these boundary conditions. 
Since  \cite{mcsalt}  the linearization
of $\mathscr D$  always falls under this heading,  we therefore have 
 the following:

\begin{prop}
For any integer $k>1$ and any real number $\alpha\in (0,1)$, 
$${\mathscr D}: C^{k+1,\alpha} (N, E)_\varepsilon \to C^{k,\alpha} (\Lambda^{0,1}\otimes N)$$
is a Fredholm map of Banach manifolds, and ${\mathcal M}={\mathscr D}^{-1}(0)$ exactly parameterizes
the holomorphic curves $(\Sigma^\prime , \partial \Sigma^\prime ) \subset (Z,P)$
which are sufficiently close to $\Sigma$.
\end{prop}

The implicit function theorem \cite{schwartz} thus implies 
an analog of Theorem \ref{koda}: 

\begin{thm}
 \label{crux}
 Let $\Sigma$ be a compact holomorphic curve-with-boundary 
 in a complex manifold $Z$, and suppose that the 
 boundary  $\partial\Sigma$ of this curve lies on a maximal
 totally real $C^5$ submanifold $P\subset Z$. If the double $\mathbb X$
 of $\Sigma$ satisfies  $H^1({\mathbb X}, {\mathcal O}({\mathscr N}))=0$, then the 
 moduli space ${\mathcal M}$ of 
 nearby holomorphic curves $(\Sigma^\prime , \partial \Sigma^\prime ) \subset (Z,P)$
is a  manifold. Moreover, the   tangent space of this manifold at the base-point $x$ 
representing $(\Sigma  , \partial \Sigma )$ 
is canonically given by 
 $T_x{\mathcal M}=H^0_\varrho ({\mathbb X}, {\mathcal O}({\mathscr N}))$.
\end{thm}

Indeed, the entire point of Lemma \ref{lift} is that $\Sigma$ 
is a regular point of 
$\mathscr D$ iff $H^1({\mathbb X}, {\mathcal O}({\mathscr N}))=0$. When
this happens, one then says that the holomorphic curve $(\Sigma, \partial \Sigma )\subset 
(Z,P)$ is {\sl Fredholm regular}. 

If we now wish to consider a $1$-parameter family of deformations of 
 either $P$ or of the complex structure of $Z$, we may do so by
simply multiplying both $C^{k+1, \alpha}(\Sigma ; N, E)$ and
$C^{k, \alpha} (\Sigma; \Lambda^{0.1}\otimes N)$ by an interval,
and augmenting the parameterized form of $\mathscr D$ with the identity
map on the second factor. 
This new map is still Fredholm, and has exactly the same
kernel and cokernel at the 
origin. The implicit function theorem therefore also yields the
following  analog of Theorem \ref{coda}:

\begin{thm}\label{redux}
Suppose $(\Sigma  , \partial \Sigma ) \subset (Z,P)$,  as above. 
 If the double $\mathbb X$
 of $\Sigma$ satisfies  $H^1({\mathbb X}, {\mathcal O}({\mathscr N}))=0$, 
then any   small deformation $(Z^\prime , P^\prime )$ of $(Z,P)$  contains
an   $h^0({\mathbb X}, {\mathcal O}({\mathscr N}))$-dimensional  family 
of holomorphic curves-with-boundary obtained by  deforming $\Sigma$. 
\end{thm}

In order to make good use of this result, one must be able to 
concretely identify both ${\mathbb X}$ and ${\mathscr N}\to {\mathbb X}$.
However, we shall now see that when $P$ is the fixed point set of 
an anti-holomorphic involution of $Z$, ${\mathbb X}$ can be identified
with a holomorphic curve immersed in $Z$, and ${\mathscr N}$
then precisely becomes the  normal bundle of ${\mathbb X}$, in the usual sense.

\section{Plane Curves}
\label{plan} 

The general theory developed in \S \ref{reflex} turns out to 
have some rather surprising consequences  for algebraic curves in the projective plane. 
 
 Let $X\subset \CP_2$ be a  complex algebraic curve 
 $$P(z_1, z_2, z_3) =0$$
 where $P$ is a homogeneous polynomial of degree ${\mathbf d}$ with {\em real} coefficients. 
 Assume that  $dP\neq 0$ along $X$, so that $X$ is smooth. Also assume
 that the real locus $C=X\cap \RR {\mathbb P}^2$ is non-empty. Let $\sigma: X\to X$
 be the anti-holomorphic involution induced by complex conjugation 
 in $\CP_2$. 
 It is not hard to see, by an elementary covering-space argument, 
  that $X-C$ has either one or two components,
 depending on whether the connected surface-with-boundary $X/\sigma$ is
 non-orientable or orientable, respectively. Both possibilities really do occur. For
 example, a real cubic can have either one or two components, and in this
 ${\mathbf d}=3$ 
 case one can check  that $X-C$ has the same number of components as the real locus.

 Let us first consider the case in which 
 $X-C$  has
 two connected components. Let $\Sigma$ be the
 closure of one of these two components. Then the corresponding double 
 ${\mathbb X}= \Sigma \cup \overline{\Sigma}$ can 
 be identified with $X$, and $\sigma: X\to X$ 
 can be identified with $\rho : {\mathbb X}\to {\mathbb X}$, and the
 virtual normal bundle ${\mathscr N}\to {\mathbb X}$ can be identified
 with the usual normal normal bundle ${\mathcal O}({\mathbf d})$ of
 $X$. Since the canonical line bundle of $X$ is ${\mathcal O}({\mathbf d}-3)$ 
 by the adjunction formula, Serre duality tells us that 
  \begin{eqnarray*}
H^1 ({\mathbb X}, {\mathcal O}({\mathscr N})) & = & H^1 ({X}, {\mathcal O}({\mathbf d}))  \\
& = & [H^0 ({X}, \Omega^1 (-{\mathbf d}))]^* \\
 & = & [H^0 ({X},  {\mathcal O} ([{\mathbf d}-3]-{\mathbf d}))]^* \\
 & = &  [H^0 ({X},  {\mathcal O} (-3))]^*  \\
 & = & 0 
\end{eqnarray*}
 since ${\mathcal O} (-3)$ has negative degree. The hypotheses of
Theorems \ref{crux} and \ref{redux} are therefore fulfilled. 

In fact, Theorem \ref{crux} tells us essentially nothing  new in this case, because the 
predicted family of dimension
$$
h^0({\mathbb X}, {\mathcal O}({\mathscr N}))= \frac{{\mathbf d}({\mathbf d}+3)}{2}= 
{{{\mathbf d}+2}\choose 2} - 1
$$
simply arises by varying the real coefficients of the homogeneous  polynomial $P$. 
However, the prediction made by Theorem \ref{redux} is, by contrast, rather
surprising. If we wiggle the embedding of $P=\RP^2$ in $Z= \CP_2$, 
 leftover halves of algebraic curves continue to cling to it, and 
the collections of ovals in $P^\prime \approx \RP^2$ which are their
boundaries give us some strange sort of deformation of the  algebraic
geometry of the real projective plane. 

When $X-C$ is connected, the  story  is basically similar, 
although   a few modest changes are necessary. 
 In this case, we instead  take 
$\Sigma$ to be a surface-with-boundary 
 diffeomorphic to $X$ minus an annular neighborhood of $C$,
 obtained from $X$ 
by  formally replacing $C$ with two disjoint copies of itself.
The double ${\mathbb X}$ of $\Sigma$  then becomes the double cover 
$\pi: \tilde{X}\to X$ given by the element of $H^1(X, \ZZ_2)$  
 Poincar\'e dual to $[C]\in H_1(X, \ZZ_2)$.  In this case, 
 we have ${\mathscr N} = \pi^* {\mathcal O}({\mathbf d})$, while
 the canonical line bundle of ${\mathbb X}= \tilde{X}$ is 
 $\pi^*{\mathcal O}({\mathbf d}-3)$, so Serre duality  tells us that 
 $$
 H^1 ({\mathbb X}, {\mathcal O}({\mathscr N}))= 
  [H^0 (\tilde{X},  \pi^*{\mathcal O} (-3))]^*=0
 $$
 because the degree of 
 $ \pi^*{\mathcal O} (-3)$ is once again negative. Thus the 
 hypotheses of
Theorems \ref{crux} and \ref{redux} are once again fulfilled. 
Here, even Theorem \ref{crux} predicts something interesting,
as 
$$
h^0({\mathbb X}, {\mathcal O}({\mathscr N}))= {\mathbf d}({\mathbf d}+3)$$
is twice as large as before, so we get deformations of real algebraic curves
which do not simply arise from varying the coefficients of $P$, but instead
arise from real deformations of the map $\tilde{X}\to \CP_2$. The observed doubling
of parameters is   analogous to what happens in the previous case if we  simultaneously
keep track of deformations
$\Sigma$ and $\overline{\Sigma}$, without requiring that their boundaries
match up in any way. Indeed, this point makes it obvious that when we
consider the  deformed
analogs of real algebraic geometry arising from the  replacement of  $P= \RP^2$ with
some nearby totally real submanifold $P^\prime\subset \CP_2$, we
must remember that the entire  story has to do with {\em oriented}
curves in $\RP^2$. In this context,  an algebraic curve with one orientation 
must  be viewed as a completely different object from its orientation-reversed twin.
 After deformation,  oppositely 
 oriented versions of any given algebraic curve will generally go their own separate ways! 
 
 Of course, the methods described here can of also be used to study real projective spaces curves
 of higher codimension , but the hypotheses of 
 Theorems \ref{crux} and \ref{redux} are  unfortunately no longer hold in general, even 
 for   complete intersections. Details are left
 to the interested reader. 

\section{Partial Indices of Holomorphic Disks}

 Many interesting applications of Theorem \ref{crux} 
already arise when   $\Sigma$ is  a disk. Under these
 circumstances,  the double of $\Sigma$ is 
 ${\mathbb X}= \CP_1$, so the Grothendieck splitting theorem \cite{oss}
 guarantees that
 $${\mathscr N}= {\mathcal O} (j_1) \oplus \cdots \oplus  {\mathcal O} (j_m).$$
If we want to consider holomorphic maps $f: (D,S^1)\to (Y,M)$ which are
not embeddings,  the theory continues to work extremely well
provided we take $Z= Y\times \CC$, $P=M\times S^1$, and let $\Sigma$
be the graph of $f$. 

In this setting, the 
numbers $j_1, \ldots , j_m$ are called the {\em partial indices}
of the disk \cite{glob}. Notice that Theorem \ref{crux} tells us \cite{ohrk} that 
such a disk is Fredholm regular iff
all the partial indices are $\geq -1$. 
The sum of the partial indices  is called the {\em Maslov index},
and corresponds to the Chern class or {\em degree} of the double
of the normal bundle.  Notice that the Maslov index is a topological 
invariant, whereas the partial indices  depend quite sensitively 
on the complex-analytic structure.

If an embedded  holomorphic  disk is real-analytic up to its boundary,
it is necessarily contained in a complex coordinate domain.
Thus, for many purposes it is quite 
sufficient to thoroughly understand  the special
case of $Z=\CC^n$. In this setting, however, the problem is amenable
to a more elementary treatment, using Fourier series or the Riemann-Hilbert
transform. In fact, many of the phenomena under discussion  
were originally discovered  \cite{forst,glob,ohrk} 
from this perspective, making them seem at  first sight  
to be completely unrelated to  
Kodaira's results. It is hoped that the present  article may play a useful r\^ole in bringing
together these  disparate strands of thought.

\section{Twistor Geometry}
 
I will now quickly describe some   twistor correspondences
involving holomorphic disks. 

Let's first return to the setting of \S \ref{plan}, and reconsider
the pair consisting of the complex manifold $Z=\CP_2$ and the totally real submanifold
$P= \RP^2\subset \CP_2$. If $C\subset \RP^2$ is any real projective line $\RP^1$,
the corresponding complex projective line $X\cong \CP_1$ is divided
into two hemispheres by $C$, and  we may single out one of these
holomorphic disks $(D, \partial D) \subset (\CP_2 , \RP^2)$ by choosing
an orientation of $C$. Thus the space $M$ of all these disks may be 
identified with the Grassmannian of oriented projective lines
in $\RP^2$, or in other words the Grassmannian
$\widetilde{Gr}_2(\RR^3)$ of 
oriented $2$-planes in $\RR^3$. This is of course just a fancy way of 
saying $S^2$, and we can clarify this  point by observing
 that each of the disks in question meets the conic 
$$z_1^2+z_2^2+z_3^2=0$$
in a unique point. This conic is of course diffeomorphic to $S^2$, and  provides a serviceable 
 model for $M$. Now, as a special case  of  
 the discussion in \S \ref{plan}, the family $M$ of holomorphic disks is stable 
 under deformations of $P$. Thus, if we wiggle the embedding
 $\RP^2\hookrightarrow \CP_2$ to produce a nearby totally real submanifold
 $P^\prime \subset \CP_2$, there is an associated $S^2$-family of holomorphic
 disks $D^\prime$ with boundaries on $P^\prime$. Let $M^\prime \approx S^2$
 be the moduli space of these disks. 
 Then $M^\prime$ contains a tautological family of closed curves.
  Indeed,  for each $y\in P^\prime$, one can 
 consider the set $L_y\subset M$ consisting of all of the
 holomorphic disks $D^\prime$ 
 passing through $y$.  
Remarkably,  the $L_y$ turn out to be exactly the 
 {\em unparameterized geodesics} of an affine connection 
$\nabla$ on $M^\prime$.  Moreover, this construction can be shown to give
rise to every connection on a compact surface for which every 
geodesic is a simple closed curve \cite{lmzoll}. 
The special case in which  $\nabla$ is the Levi-Civita connection of 
some Riemannian metric $g$ can also be thoroughly analyzed from this point of view, 
leading  to an entirely new 
understanding \cite{lmzoll} of 
the classical theory of  Zoll surfaces \cite{beszoll}.

 Next, let us consider what happens if we instead
 take $Z=\CP_3$ and $P= \RP^3\subset \CP_3$. 
 Again, every real projective line bounds two
 holomorphic disks, and the moduli space of these
 disks is now  the Grassmannian
$\widetilde{Gr}_2(\RR^4)$ of 
oriented $2$-planes in $\RR^4$. Each such disk 
meets the quadric 
$$z_1^2+z_2^2+z_3^2+z_4^2=0$$
in a unique point, so this oriented Grassmannian 
can be identified, if we like, with a complex
$2$-quadric $Q_2\approx S^2\times S^2$. 
Because  the double ${\mathbb X}$ of any
such disk is a projective line  $\CP_1$ with normal bundle
${\mathscr N}\cong {\mathcal O}(1) \oplus {\mathcal O}(1)$, 
these disks  all satisfy $H^1({\mathbb X}, {\mathcal O}({\mathscr N}))=0$. 
Theorem \ref{redux} thus predicts that if we perturb 
$\RP^3\hookrightarrow \CP_3$  to obtain 
a nearby totally real submanifold $P^\prime \subset \CP_3$,
there is an analogous $(S^2\times S^2)$-family 
of holomorphic disks $D^\prime$ with boundaries in $P^\prime$.
Let $M^\prime \approx S^2\times S^2$ denote this moduli
space. Then $M^\prime$ comes equipped with a
natural family of embedded $2$-spheres $S_y\subset M^\prime$,
where, for each $y\in P^\prime$, $S_y$ consists of all  the
disks in the family passing through $y$. One can then show
that there is a  pseudo-Riemannian metric $g$ on $M$
with respect to which the $S_y$ are all null surfaces. 
This $g$ is unique up to conformal rescaling, and is {\em self-dual},
in the sense that its Weyl curvature $W$ satisfies
$\star W= W$ as a bundle-valued $2$-form. 
Moreover, every self-dual conformal metric near the 
standard one arises from this construction
\cite{lmfrei}. Notice that the geometries that arise from this construction, 
unlike the previous one, now satisfy a local curvature condition. 
On the other hand, global conditions on the periodicity of
geodesics do not have to be explicitly stipulated in this case, as they turn out
to automatically hold for any solution which is $C^2$ close to the 
standard one.

It is now only natural to  ask what happens if  we instead take 
$Z=\CP_{m+1}$ and $P= \RP^{m+1}\subset \CP_{m+1}$
for some $m\geq 3$. While this is a story which has never properly been 
set down in detail, many of the broad outlines are certainly similar to what we
have already seen. 
Every real projective line in $\RP^{m+1}$ bounds two
 holomorphic disks in $\CP_{m+1}$, and the moduli space of these
 disks is   the Grassmannian
$\widetilde{Gr}_2(\RR^{m+2})$ of 
oriented $2$-planes in $\RR^{m+2}$. Each such disk 
meets the quadric 
$$z_1^2+z_2^2+\cdots +z_{m+2}^2=0$$
in a unique point, so this oriented Grassmannian 
can be identified  with the complex
$m$-quadric $Q_m$. 
Because  the double ${\mathbb X}$ of any
such disk is a  projective line with normal bundle
${\mathscr N}\cong [{\mathcal O}(1)]^{\oplus m}$, 
these disks are all Fredholm regular, and 
Theorem \ref{redux} again tells us that each 
 totally real submanifold $P^\prime$
 near the standard 
$\RP^{m+1}\subset \CP_{m+1}$  
has an associated  family 
of holomorphic disks $D^\prime$ with boundaries in $P^\prime$.
The  moduli
space $M^\prime \approx Q_m$ 
of these disks 
contains a
tautological  family of embedded $m$-spheres $S_y\subset M^\prime$,
 $y\in P^\prime$, given by  sub-families of those
disks   passing through any given $y$. This time, however, 
the associated geometry of 
$M^\prime$ is, in the terminology of  \cite{paraconf},
 exactly  a right-flat $(m,2)$-paraconformal
structure; when $m$ is even, this sort of structure may 
be thought of as a Wick-rotated version of a 
  quaternionic structure   \cite{salquat}.
Given any volume form on $M^\prime\approx Q_m$,
there is a unique torsion-free 
affine connection of holonomy $\subset [SL(2, \RR) \times SL(m, \RR)]/\ZZ_2$
 which is compatible with the paraconformal structure
and the volume form;  the submanifolds $S_y$ are then
totally geodesic with respect to this connection. 
The arguments in \cite{lmfrei} strongly indicate
 that every such structure on $Q_m$ sufficiently near  the  
standard one should arise from this construction. In particular,  
the general such structure 
on $Q_m$   should   depend on $(m+1)$ real functions
of $(m+1)$ real variables. Details are 
 left to the interested reader. 

These are but a few simple examples of the manner 
 in which  moduli of holomorphic  curves-with-boundary can naturally
give rise to geometrically rich twistor correspondences.  Of course, 
we have in each case simply taken $Z$ to be $\CP_n$, and taken $\Sigma$ to
be half a projective line. There are certainly many, many more correspondences
of this same flavor, just waiting    to be developed. For example, 
 by  taking $\Sigma$ to be a disk in a different  rational complex surface,
one would  encounter versions of  the Hitchin  correspondence for  Einstein-Weyl spaces
\cite{hitproj} 
or Bryant's connections with  exotic affine holonomy \cite{bryex}. 
But what if we take $\Sigma$ to be something other than a disk?
Such  moduli spaces must  certainly carry fascinating geometries 
whose secrets are simply  waiting  to be unlocked. I can only hope that 
some interested reader will take up the challenge, and try to 
chart a bit  of this  {\em terra incognita}.

  \end{document}